\documentclass{amsart}

\usepackage{amssymb}
\usepackage{amscd}

  \theoremstyle{plain}
  \newtheorem*{thm*}{Theorem}
  \theoremstyle{plain}
  \newtheorem{thm}{Theorem}[section]
  \theoremstyle{definition}
  \newtheorem{defn}[thm]{Definition}
  \theoremstyle{plain}
  
  \theoremstyle{plain}
  
  \theoremstyle{plain}
  \newtheorem{cor}[thm]{Corollary}
  \theoremstyle{remark}
  
  \theoremstyle{remark}
  \newtheorem*{acknowledgement*}{Acknowledgement}

\begin{document}

\title[Funk functions for Finsler spaces]{Non-existence of Funk functions for Finsler spaces of
  non-vanishing scalar flag curvature.}

\author[Bucataru]{Ioan Bucataru} \address{Ioan Bucataru, Faculty of
Mathematics, Alexandru Ioan Cuza University \\ Ia\c si, 
Romania} \urladdr{http://www.math.uaic.ro/$\sim$bucataru/}

\author[Muzsnay]{Zolt\'an Muzsnay},
\address{Institute of Mathematics, University of Debrecen, Debrecen, Hungary}
\urladdr{http://math.unideb.hu/$\sim$muzsnay-zoltan}

\date{\today}

\begin{abstract}
In his book "Differential Geometry of Spray and Finsler spaces", page 177, Zhongmin Shen asks ``wether or not there
always exist non-trivial Funk functions on a spray space''. In this
note, we will prove that the answer is negative for the geodesic spray of a
finslerian function of non-vanishing scalar flag curvature. 
\end{abstract}



\maketitle

\section{Introduction}\label{introduction}

A Funk function is a projective deformation of a homogeneous system of second
order ordinary differential equations (SODE) that preserves the
curvature (Jacobi endomorphism) of the given SODE. 

If we start with a flat SODE, there are many corresponding Funk
functions and each of them will projectively deform the given SODE
into a Finsler metrizable one, \cite[Theorem 7.1]{GM00},
\cite[Theorem 10.3.5]{Shen01}, solution to Hilbert's fourth problem. 

The existence of non-trivial Funk functions, for the general case of a
non-flat SODE, is an open problem, as mentioned by Zhongim Shen in \cite[page 177]{Shen01}. 
We will prove in Theorem \ref{thm1} that there are non-trivial Funk
functions for the geodesic equations of a Finsler function of
non-vanishing scalar flag curvature.

\section{Funk functions and projective deformations of Finsler spaces}

We consider $M$ a smooth $n$-dimensional manifold, $TM$ its tangent
bundle and $T_0M=TM\setminus \{0\}$ the tangent bundle with the zero
section removed. Local coordinates on $M$ are denoted by $(x^i)$, while induced
local coordinates on $TM$ and $T_0M$ are denoted by $(x^i, y^i)$, for $i\in\{1,...,
n\}$. 

Since most of the geometric objects that we will use in this paper are
homogeneous with respect to the fibre coordinates on $TM$, we will
assume them to be defined on $T_0M$. The homogeneity is characterised
with the help of the Liouville (or dilation) vector field ${\mathbb
  C}=y^i{\partial}/{\partial y^i}$.

The geometric framework that we will use in this paper is based on the
Fr\"olicher-Nijenhuis formalism, \cite[Chapter 2]{GM00}, \cite{Szilasi03}. Within this formalism
one can associate to each vector valued form $L$,
of degree $\ell$, two derivations: $i_L$, of degree $(\ell-1)$ and
$d_L$, of degrees $\ell$. For two vector valued forms $K$ and $L$, of
degree $k$ and $\ell$, we consider the Fr\"olicher-Nijenhuis bracket
$[K,L]$, which is a vector valued $(k+\ell)$-form, uniquely determined
by $d_{[K,L]}=d_Kd_L-(-1)^{k\ell}d_Ld_K$.
For various commutation formulae regarding the Fr\"olicher-Nijenhuis
theory, we refer to \cite[Appendix A]{GM00}. 
On $TM$ there is a canonical vector valued
$1$-form, the tangent endomorphism $J=dx^i\otimes {\partial}/{\partial
  y^i}$ that induces two
derivations $i_J$ and $d_J$.

A system of second order ordinary differential equations on $M$,
\begin{eqnarray}
\frac{d^2x^i}{dt^2} + 2G^i\left(x,\frac{dx}{dt}\right) =
  0, \label{sode} \end{eqnarray}
can be identified with a vector field
$S \in {\mathfrak X}(TM)$,  $S=y^i{\partial}/{\partial x^i} - 2G^i(x,y){\partial}/{\partial
  y^i}$, which satisfies $JS={\mathbb C}$. Such a vector field is
called a semispray. If additionally, $S\in
{\mathfrak X}(T_0M)$ and satisfies the homogeneity condition $[{\mathbb
    C}, S]=S$ we say that $S$ is a \emph{spray}.

For a spray $S$, we can associate a nonlinear (Ehresman)
connection with curvature tensors, \cite[\S
3.1, \S 3.2]{GM00}. We denote by $h$ the horizontal
projector, $R$ the curvature of the connection and
$\Phi$ the Jacobi endomorphism, given by 
\begin{eqnarray*} 
h=\frac{1}{2}\left( Id - [S,J]\right), \quad R=\frac{1}{2}[h,h], \quad
  \Phi = \left( Id - h \right) \circ [S,h]. \end{eqnarray*}
The two curvature tensors $R$ and $\Phi$ and the induced derivations are
related by:
\begin{eqnarray*} \Phi=i_SR, \quad [J, \Phi]=3R, \quad
3d_R=d_{[J, \Phi]} = d_Jd_{\Phi} + d_{\Phi}d_J. \end{eqnarray*}
A spray $S\in {\mathfrak X}(T_0M)$ is said to be \emph{isotropic} if
the induced Jacobi endomorphism takes the particular form 
\begin{eqnarray}
\Phi = \rho J - \alpha \otimes {\mathbb
  C}. \label{Phi_iso} \end{eqnarray} 
The function $\rho \in C^{\infty}(T_0M)$ is called the \emph{Ricci scalar}
and $\alpha=\alpha_i dx^i$ is a semi-basic $1$-form (with respect to
the canonical projection of the tangent bundle).

\begin{defn}
A \emph{Finsler function} is a continuous non-negative function $F: TM\to {\mathbb R}$ that
satisfies the following conditions:
\begin{itemize} \label{def_Finsler}
\item[i)] $F$ is smooth on $T_0M$ and $F(x,y)=0$ if and only if $y=0$;
\item[ii)] $F$ is positively homogeneous of order $1$ in the fiber
  coordinates;
\item[iii)] the $2$-form $dd_JF^2$ is a symplectic form on $T_0M$.
\end{itemize} 
\end{defn}
A spray $S\in {\mathfrak X}(T_0M)$ is said to be \emph{Finsler metrizable} if
there exists a Finsler function $F$ that satisfies
\begin{eqnarray}
i_Sdd_JF^2=-dF^2. \label{geodspray1}
\end{eqnarray}
For a given Finsler function $F$, equation (\ref{geodspray1}) uniquely determine a spray $S$, called the \emph{geodesic spray} of the Finsler
function $F$. The geometric structures associated to a Finsler
function are those corresponding to its geodesic spray. Equation
(\ref{geodspray1}) is equivalent to $ d_h F^2 = 0.$ 

A Finsler function $F$ has \emph{scalar flag curvature} $\kappa\in C^{\infty}(T_0M)$ if
the Jacobi endomorphism is given by  
\begin{eqnarray}
\Phi=\kappa\left (F^2 J - Fd_JF \otimes {\mathbb
  C}\right). \label{sfc}
\end{eqnarray}

An orientation-preserving reparameterization of the system
(\ref{sode}) gives rise to a new system and therefore a new spray $\widetilde{S}=S-2P{\mathbb C}$, \cite[Chapter
12]{Shen01}. These two sprays $S$ and $\widetilde{S}$ are called
\emph{projectively related}. The function $P$ is $1$-homogeneous and
it is called the projective deformation of the spray $S$. 

Under a projective deformation $\widetilde{S}=S-2P{\mathbb C}$, the
corresponding geometric setting changes, the Jacobi
endomorphisms being related by, \cite[(4.8)]{BM12}, 
\begin{eqnarray}
\widetilde{\Phi} & = & \Phi + \left(P^2 - S(P)\right) J - \left(d_J(S(P) - P^2) +
  3(Pd_JP  - d_{h}P) \right)\otimes {\mathbb C}. \label{phiphi0}
\end{eqnarray}
 From the above formula (\ref{phiphi0}) we obtain that a projective
 deformation $P$ preserves the Jacobi endomorphism if and only if it
 is a \emph{Funk function}, which means that it satisfies  
\begin{eqnarray}
d_hP=Pd_JP. \label{Funk}
\end{eqnarray}     

\section{Non-existence results for Funk functions}

We will prove now that there is an obstruction
for the Funk equation (\ref{Funk}), which is not satisfied for Finsler
functions of non-vanishing scalar flag curvature. 

\begin{thm} \label{thm1}
Consider $F$ a Finsler function $F$ of non-vanishing scalar flag curvature
$\kappa$. Then, there are no non-trivial Funk functions for the Finsler space $(M, F)$.
\end{thm}
\textit{Proof.}
Consider $S$ the geodesic spray of a Finsler function $F$, of non-vanishing
scalar flag curvature $\kappa$. It follows that its Jacobi endomorphism is given
by formula (\ref{sfc}).

We will prove the statement of the theorem by contradiction. Therefore, we
assume, that there exists a Funk function $P$, for the Finsler space
$(M, F)$. Hence, the function $P$ satisfies the Funk equation (\ref{Funk}).

If we apply the derivation
$d_J$ to both sides of (\ref{Funk}) we obtain $d_Jd_{h}P=0$. Since
$[J,h]=0$, we obtain $0=d_{[J,h]}P=d_Jd_hP+d_hd_JP$ and hence
$d_hd_JP=0$. Using this formula and 
if we apply $d_{h}$ to both sides of (\ref{Funk}) we obtain 
\begin{eqnarray}
d_{R}P=d^2_{h}P=d_{h}\left(Pd_JP\right)=d_h P\wedge d_JP + Pd_hd_JP= 0. \label{dRP}
\end{eqnarray}
Due to the homogeneity of the involved geometric structures and the
relations between the Jacobi endomorphism $\Phi$ and the curvature $R$ of the nonlinear
connection $h$, we obtain that the equation (\ref{dRP}) is
equivalent to 
\begin{eqnarray}
0=i_Sd_RP=d_{i_SR}P=d_{\Phi}P=0. \label{dPhiP}
\end{eqnarray}
Using the form (\ref{sfc}) of the Jacobi endomorphism $\Phi$, the
above equation (\ref{dPhiP}) can be written as follows
\begin{eqnarray}
\kappa F^2d_JP - \kappa F Pd_JF=0 \Longleftrightarrow
  d_J\left(\frac{P}{F}\right)=0. \label{djpf}
\end{eqnarray}
Last equation above implies that the function $P/F$ is constant along
the fibres the tangent bundle. Therefore, there exists a basic
function $a$, on the base manifold $M$, such that 
\begin{eqnarray}
P(x,y)=a(x)F(x,y), \forall (x,y) \in T_0M. \label{pafin}
\end{eqnarray}
Now, if we use the above form of the projective factor $P$ and the fact
that $d_hF=0$, the Funk equation (\ref{Funk}) can be written as follows
\begin{eqnarray}
Fda=a^2 F d_JF \Longleftrightarrow
  d\left(\frac{-1}{a}\right)=d_JF. \label{dajf}
\end{eqnarray}
The last equation (\ref{dajf}) cannot have solutions since the left
hand side is a basic $1$-form, while the right hand side is a
semi-basic (and not basic) $1$-form. \qed

We can reinterpret Theorem \ref{thm1} as follows. Consider $S$ the geodesic spray of a Finsler
function of non-vanishing scalar flag curvature. Then there is no
projective deformation of the spray $S$, by a projective factor $P$, which
preserves the curvature of the spray. 

In \cite[Theorem 3.1]{Bucataru15} it is shown that for a Finsler
function of non-vanishing scalar flag curvature, the projective class
of its geodesic spray contains exactly one spray with the given curvature
tensor. Next corollary, which is a consequence of Theorem \ref{thm1}, gives a new
motivation for \cite[Theorem 3.1]{Bucataru15}. 

\begin{cor} \label{cor1}
Consider $S$ an isotropic spray of non-vanishing Ricci scalar. Then
there is no Funk function on the spray space $(M, S)$ that will
projectively deform $S$ into a Finsler metrizable spray.
\end{cor}
\textit{Proof.}
Consider $S$ an isotropic spray of non-vanishing Ricci scalar. By
contradiction, we assume that there is a Funk function $P$, for the
spray space $(M, S)$, such that
the projectively deformed spray $\widetilde{S}=S-2P{\mathbb C}$ is
Finsler metrizable by a Finsler function $F$. Since $P$ is a Funk
function it follows that the two sprays $S$ and $\widetilde{S}$ have
the same Jacobi endomorphisms and hence $\widetilde{S}$ is isotropic
as well, with non-vanishing Ricci scalar. Since $\widetilde{S}$ is Finsler metrizable by a Finsler
function $F$, it follows that $F$ is of non-vanishing scalar flag
curvature. It follows that $-P$ is a Funk function for the Finsler
space $(M, F)$, which contradicts the result of Theorem \ref{thm1}. 
\qed

The results of Theorem \ref{thm1} and Corollary \ref{cor1} heavily
rely on the fact that the curvature of the considered spray and Finsler
function does not vanish. This implies that for the existence of the
Funk function, we have the obstructions (\ref{dRP}) and
(\ref{dPhiP}). In the flat case, these obstructions are automatically
satisfied, and it explains why in this case there are many Funk
functions and each of them leads to sprays that are Finsler metrizable, \cite[Theorem 7.1]{GM00},
\cite[Theorem 10.3.5]{Shen01}.

\section*{Acknowledgements}
This work was supported by the Bilateral Cooperation Programs: RO-HU
672/2013-2014, T\'ET-12-RO-1-2013-0022 and EU FET FP7 BIOMICS project:
CNECT-318202.

\end{document}